\documentclass{amsproc}
\usepackage{amssymb}
\usepackage{amsmath}
\usepackage{eurosym}
\usepackage{amsfonts}

            \advance \topmargin by -\headheight \advance
            \topmargin by -\headsep
            \evensidemargin
            \oddsidemargin
\newtheorem{theorem}{Theorem}

\newtheorem{corollary}[theorem]{Corollary}

\newtheorem{lemma}[theorem]{Lemma}

\newtheorem{remark}[theorem]{Remark}

\numberwithin{equation}{section}
\input{tcilatex}

\begin{document}
\title[Tetranacci and Tetranacci-Lucas Quaternions]{Tetranacci and
Tetranacci-Lucas Quaternions}
\thanks{}
\author[Y\"{u}ksel Soykan]{Y\"{u}ksel Soykan}
\maketitle

\begin{center}
\textsl{Department of Mathematics, Art and Science Faculty, }

\textsl{Zonguldak B\"{u}lent Ecevit University, 67100, Zonguldak, Turkey }
\end{center}

\textbf{Abstract.} The quaternions form a $4$-dimensional Cayley-Dickson
algebra. In this paper, we introduce the Tetranacci and Tetranacci-Lucas
quaternions. Furthermore, we present some properties of these quaternions
and derive relationships between them.



\section{Introduction}

Tetranacci sequence $\{M_{n}\}_{n\geq 0}$ and Tetranacci-Lucas sequence $%
\{R_{n}\}_{n\geq 0}$ are defined by the fourth-order recurrence relations%
\begin{equation}
M_{n}=M_{n-1}+M_{n-2}+M_{n-3}+M_{n-4},\text{ \ \ \ \ }%
M_{0}=0,M_{1}=1,M_{2}=1,M_{3}=2  \label{equati:fvcvxghsbnz}
\end{equation}%
and 
\begin{equation}
R_{n}=R_{n-1}+R_{n-2}+R_{n-3}+R_{n-4},\text{ \ \ \ \ }%
R_{0}=4,R_{1}=1,R_{2}=3,R_{3}=7  \label{equati:pazertvbcunsmn}
\end{equation}%
respectively. $M_{n}$ is the sequence $A000078$ in [\ref{bib:sloane}] and $%
R_{n}$ is the sequence $A073817$ in [\ref{bib:sloane}]. This sequence has
been studied by many authors and more detail can be found in the extensive
literature dedicated to these sequences, see for example [\ref%
{bib:hathiwala2017}], [\ref{melham1999}], [\ref{natividad2013}], [\ref%
{bib:singh2014}], [\ref{waddill1967}], [\ref{waddill1992}].

The sequences $\{M_{n}\}_{n\geq 0}$ and $\{R_{n}\}_{n\geq 0}$ can be
extended to negative subscripts by defining 
\begin{equation*}
M_{-n}=-M_{-(n-1)}-M_{-(n-2)}-M_{-(n-3)}+M_{-(n-4)}
\end{equation*}%
and%
\begin{equation*}
R_{-n}=-R_{-(n-1)}-R_{-(n-2)}-R_{-(n-3)}+R_{-(n-4)}
\end{equation*}%
for $n=1,2,3,...$ respectively. Therefore, recurrences (\ref%
{equati:fvcvxghsbnz}) and (\ref{equati:pazertvbcunsmn}) hold for all integer 
$n.$

We can write (\ref{equati:fvcvxghsbnz}) as $%
M_{n-1}=M_{n-2}+M_{n-3}+M_{n-4}+M_{n-5}.$ Substracting this from (\ref%
{equati:fvcvxghsbnz}), we see that Tetranacci numbers also satisfy the
following useful alternative linear recurrence relation for $n\geq 5$:%
\begin{equation}
M_{n}=2M_{n-1}-M_{n-5}.  \label{equation:sdfghbcxvzfh}
\end{equation}

Extension of the definition of $M_{n}$ to negative subscripts can be proved
by writing the recurrence relation (\ref{equation:sdfghbcxvzfh}) as%
\begin{equation}
M_{-n}=2M_{-n+5}-M_{-n+6}.  \label{equation:yuhbncvxgfsrta}
\end{equation}

Similarly, we have%
\begin{eqnarray}
R_{n} &=&2R_{n-1}-R_{n-5}, \\
R_{-n} &=&2R_{-n+5}-R_{-n+6}.
\end{eqnarray}

The following Table 1 presents the first few values of the Tetranacci and
Tetranacci-Lucas numbers with positive and negative subscripts:

\ \ \ \ 

Table 1. Tetranacci and Tetranacci-Lucas Numbers with non-negative and
negative indices

\begin{tabular}{cccccccccccccccc}
\hline
$n$ & $0$ & $1$ & $2$ & $3$ & $4$ & $5$ & $6$ & $7$ & $8$ & $9$ & $10$ & $11$
& $12$ & $13$ & $...$ \\ \hline
$M_{n}$ & $0$ & $1$ & $1$ & $2$ & $4$ & $8$ & $15$ & $29$ & $56$ & $108$ & $%
208$ & $401$ & $773$ & $1490$ & $...$ \\ 
$M_{-n}$ & $0$ & $0$ & $0$ & $1$ & $-1$ & $0$ & $0$ & $2$ & $-3$ & $1$ & $0$
& $4$ & $-8$ & $5$ & $...$ \\ 
$R_{n}$ & $4$ & $1$ & $3$ & $7$ & $15$ & $26$ & $51$ & $99$ & $191$ & $367$
& $708$ & $1365$ & $2631$ & $5071$ & $...$ \\ 
$R_{-n}$ & $4$ & $-1$ & $-1$ & $-1$ & $7$ & $-6$ & $-1$ & $-1$ & $15$ & $-19$
& $4$ & $-1$ & $31$ & $-53$ & $...$ \\ \hline
\end{tabular}

\ \ \ \ \ \ \ \ \ 

It is well known that for all integers $n,$ usual Tetranacci and
Tetranacci-Lucas numbers can be expressed using Binet's formulas%
\begin{equation*}
M_{n}=\frac{\alpha ^{n+2}}{(\alpha -\beta )(\alpha -\gamma )(\alpha -\delta )%
}+\frac{\beta ^{n+2}}{(\beta -\alpha )(\beta -\gamma )(\beta -\delta )}+%
\frac{\gamma ^{n+2}}{(\gamma -\alpha )(\gamma -\beta )(\gamma -\delta )}+%
\frac{\delta ^{n+2}}{(\delta -\alpha )(\delta -\beta )(\delta -\gamma )}
\end{equation*}%
(see for example [\ref{bib:zaveri2015}] or [\ref{bib:hathiwala2017}])

or%
\begin{equation}
M_{n}=\frac{\alpha -1}{5\alpha -8}\alpha ^{n-1}+\frac{\beta -1}{5\beta -8}%
\beta ^{n-1}+\frac{\gamma -1}{5\gamma -8}\gamma ^{n-1}++\frac{\delta -1}{%
5\delta -8}\delta ^{n-1}  \label{equation:srtynsbfdtsra}
\end{equation}%
(see for example [\ref{dresden2014}]) and%
\begin{equation*}
R_{n}=\alpha ^{n}+\beta ^{n}+\gamma ^{n}+\delta ^{n}
\end{equation*}%
respectively, where $\alpha ,\beta ,\gamma $ and $\delta $ are the roots of
the cubic equation $x^{4}-x^{3}-x^{2}-x-1=0.$ Moreover,%
\begin{eqnarray*}
\alpha &=&\frac{1}{4}+\frac{1}{2}\omega +\frac{1}{2}\sqrt{\frac{11}{4}%
-\omega ^{2}+\frac{13}{4}\omega ^{-1}}, \\
\beta &=&\frac{1}{4}+\frac{1}{2}\omega -\frac{1}{2}\sqrt{\frac{11}{4}-\omega
^{2}+\frac{13}{4}\omega ^{-1}}, \\
\gamma &=&\frac{1}{4}-\frac{1}{2}\omega +\frac{1}{2}\sqrt{\frac{11}{4}%
-\omega ^{2}-\frac{13}{4}\omega ^{-1}}, \\
\delta &=&\frac{1}{4}-\frac{1}{2}\omega -\frac{1}{2}\sqrt{\frac{11}{4}%
-\omega ^{2}-\frac{13}{4}\omega ^{-1}},
\end{eqnarray*}%
where%
\begin{equation*}
\omega =\sqrt{\frac{11}{12}+\left( \frac{-65}{54}+\sqrt{\frac{563}{108}}%
\right) ^{1/3}+\left( \frac{-65}{54}-\sqrt{\frac{563}{108}}\right) ^{1/3}}.
\end{equation*}%
Note that we have the following identities: 
\begin{eqnarray*}
\alpha +\beta +\gamma +\delta &=&1, \\
\alpha \beta +\alpha \gamma +\alpha \delta +\beta \gamma +\beta \delta
+\gamma \delta &=&-1, \\
\alpha \beta \gamma +\alpha \beta \delta +\alpha \gamma \allowbreak \delta
+\beta \gamma \delta &=&1, \\
\alpha \beta \gamma \delta &=&-1.
\end{eqnarray*}

Note that the Binet form of a sequence satisfying (\ref{equati:fvcvxghsbnz}%
)\ and (\ref{equati:pazertvbcunsmn}) for non-negative integers is valid for
all integers $n.$ This result of Howard and Saidak [\ref{bib:howard2010}] is
even true in the case of higher-order recurrence relations as the following
theorem shows.

\begin{theorem}[{[\protect\ref{bib:howard2010}]}]
\label{theorem:fvgxdfsxczsaer}Let $\{w_{n}\}$ be a sequence such that 
\begin{equation*}
\{w_{n}\}=a_{1}w_{n-1}+a_{2}w_{n-2}+...+a_{k}w_{n-k}
\end{equation*}%
for all integers $n,$ with arbitrary initial conditions $%
w_{0},w_{1},...,w_{k-1}.$ Assume that each $a_{i}$ and the initial
conditions are complex numbers. Write%
\begin{eqnarray}
f(x) &=&x^{k}-a_{1}x^{k-1}-a_{2}x^{k-2}-...-a_{k-1}x-a_{k}
\label{equation:mnbvyuhgoewapvbc} \\
&=&(x-\alpha _{1})^{d_{1}}(x-\alpha _{2})^{d_{2}}...(x-\alpha _{h})^{d_{h}} 
\notag
\end{eqnarray}%
with $d_{1}+d_{2}+...+d_{h}=k,$ and $\alpha _{1},\alpha _{2},...,\alpha _{k}$
distinct. Then

\begin{description}
\item[(a)] For all $n,$%
\begin{equation}
w_{n}=\sum_{m=1}^{k}N(n,m)(\alpha _{m})^{n}
\label{equation:yuhnbvdfscxzratqw}
\end{equation}%
where%
\begin{equation*}
N(n,m)=A_{1}^{(m)}+A_{2}^{(m)}n+...+A_{r_{m}}^{(m)}n^{r_{m}-1}=%
\sum_{u=0}^{r_{m}-1}A_{u+1}^{(m)}n^{u}
\end{equation*}%
with each $A_{i}^{(m)}$ a constant determined by the initial conditions for $%
\{w_{n}\}$. Here, equation (\ref{equation:yuhnbvdfscxzratqw}) is called the
Binet form (or Binet formula) for $\{w_{n}\}.$ We assume that $f(0)\neq 0$
so that $\{w_{n}\}$ can be extended to negative integers $n.$

If the zeros of (\ref{equation:mnbvyuhgoewapvbc}) are distinct, as they are
in our examples, then%
\begin{equation*}
w_{n}=A_{1}(\alpha _{1})^{n}+A_{2}(\alpha _{2})^{n}+...+A_{k}(\alpha
_{k})^{n}.
\end{equation*}

\item[(b)] The Binet form for $\{w_{n}\}$ is valid for all integers $n.$
\end{description}
\end{theorem}

The generating functions for the Tetranacci sequence $\{M_{n}\}_{n\geq 0}$
and Tetranacci-Lucas sequence $\{R_{n}\}_{n\geq 0}$ are 
\begin{equation*}
\sum_{n=0}^{\infty }M_{n}x^{n}=\frac{x}{1-x-x^{2}-x^{3}-x^{4}}\text{ \ and \ 
}\sum_{n=0}^{\infty }R_{n}x^{n}=\frac{4-3x-2x^{2}-x^{3}}{%
1-x-x^{2}-x^{3}-x^{4}},
\end{equation*}%
respectively.

In this paper, we define Tetranacci and Tetranacci-Lucas quaternions in the
next section and give some properties of them. Before giving their
definition, we present some information on quaternions.

Quaternions were invented by Irish mathematician W. R. Hamilton (1805-1865)
as an extension to the complex numbers. Most mathematicians have heard the
story of how Hamilton invented the quaternions. The 16th of October 1843 was
a momentous day in the history of mathematics and in particular a major
turning point in the subject of algebra. On that day William Rowan Hamilton
had a brain wave and came up with the idea of the quaternions. He carved the
multiplication formulae with his knife into the stone of the Brougham Bridge
(nowadays known as Broomebridge) in Dublin,%
\begin{equation*}
i^{2}=j^{2}=k^{2}=ijk=-1.
\end{equation*}%
One reason this story is so well-known is that Hamilton spent the rest of
his life obsessed with the quaternions and their applications to geometry.
The story of this discovery has been translated into many different
languages. For this story and for a full biography of Hamilton, we refer the
work of Hankins [\ref{hankins1980}].

After the middle of the 20th century, the practical use of quaternions has
been discovered in comparison with other methods and there has been an
increasing interest in algebra problems on quaternion field since many
algebra problems on quaternion field were encountered in some applied and
pure science such as the quantum physics, computer science, analysis and
differential geometry.

A quaternion is a hyper-complex number and is defined by

\begin{equation*}
q=a_{0}+ia_{1}+ja_{2}+ka_{3}=(a_{0},a_{1},a_{2},a_{3})
\end{equation*}%
where $a_{0},a_{1},a_{2}\ $and $a_{3}$ are real numbers or scalers and $%
1,i,j,k$ are the standard orthonormal basis in $\mathbb{R}^{4}$. The set of
all quaternions are denoted by $\mathbb{H}$. Note that we can write%
\begin{equation*}
q=a_{0}+p
\end{equation*}%
where $p=ia_{1}+ja_{2}+ka_{3}$. $a_{0}$ and $p$ are called the scalar part
and the vector part of the quaternion $q,$ respectively. The $%
a_{0},a_{1},a_{2},a_{3}$ are called the components of the quaternion $q$.

Addition of quaternions is defined as componentwise and the quaternion
multiplication is defined as follows:%
\begin{equation}
i^{2}=j^{2}=k^{2}=ijk=-1.  \label{equati:nbvcdftyrbvcx}
\end{equation}%
Note that from (\ref{equati:nbvcdftyrbvcx}), we have%
\begin{equation}
ij=k=-ji,\text{ }jk=i=-kj,\text{ }ki=j=-ik.  \label{equation:yunbasdfxczraew}
\end{equation}%
So, multiplication on $\mathbb{H}$ is not commutative. The identities in (%
\ref{equati:nbvcdftyrbvcx}) and (\ref{equation:yunbasdfxczraew}), sometimes
are known as Hamilton's rules. Quaternions have the following multiplication
Table 2:

\ \ \ \ \ \ \ \ \ \ \ 

Table 2. Multiplication Table

$%
\begin{tabular}[t]{|l|l|l|l|l|}
\hline
& $1$ & $i$ & $j$ & $k$ \\ \hline
$1$ & $1$ & $i$ & $j$ & $k$ \\ \hline
$i$ & $i$ & $-1$ & $k$ & $-j$ \\ \hline
$j$ & $j$ & $-k$ & $-1$ & $i$ \\ \hline
$k$ & $k$ & $j$ & $-i$ & $-1$ \\ \hline
\end{tabular}%
$

\ \ \ $\ $

The product of two quaternions $q=a_{0}+ia_{1}+ja_{2}+ka_{3}$ and $%
p=b_{0}+ib_{1}+jb_{2}+kb_{3}$\ is%
\begin{eqnarray*}
qp
&=&(a_{0}b_{0}-a_{1}b_{1}-a_{2}b_{2}-a_{3}b_{3})+i(a_{0}b_{1}+a_{1}b_{0}+a_{2}b_{3}-a_{3}b_{2})
\\
&&+j(a_{0}b_{2}-a_{1}b_{3}+a_{2}b_{0}+a_{3}b_{1})+k(a_{0}b_{3}+a_{1}b_{2}-a_{2}b_{1}+a_{3}b_{0}).
\end{eqnarray*}%
The conjugate of the quaternion $q$ is defined by 
\begin{equation*}
q^{\ast }=(a_{0}+ia_{1}+ja_{2}+ka_{3})^{\ast }=a_{0}-ia_{1}-ja_{2}-ka_{3}.
\end{equation*}%
For two quaternions $p,q$ we have%
\begin{equation*}
(q^{\ast })^{\ast }=q,\text{ }(p+q)^{\ast }=p^{\ast }+q^{\ast },\text{ }%
(pq)^{\ast }=q^{\ast }p^{\ast }\text{ and }(p^{\ast }q)^{\ast }=q^{\ast }p.%
\text{ }
\end{equation*}%
The norm of a quaternion $q$ is defined by%
\begin{equation*}
N(q)=\left\Vert q\right\Vert :=qq^{\ast
}=a_{0}^{2}+a_{1}^{2}+a_{2}^{2}+a_{3}^{2}.
\end{equation*}%
The norm is multiplicative:%
\begin{equation*}
N(pq)=N(p)N(q).
\end{equation*}%
Division is uniquely defined (except by zero), thus quaternions form a
division algebra. For two quaternions $p,q\in \mathbb{H}$ we have 
\begin{equation*}
(pq)^{-1}=q^{-1}p^{-1}.
\end{equation*}%
The inverse (reciprocal) of a nonzero quaternion\ $q$ is given by%
\begin{equation*}
q^{-1}=\frac{q^{\ast }}{N(q)}.
\end{equation*}

In $1898$ A. Hurwitz proved that the only real composition algebras are $%
\mathbb{R}$, $\mathbb{C}$, $\mathbb{H}$ and $\mathbb{O}$ (here $\mathbb{O}$
stands for octonion algebras). (A real composition algebra is an algebra $%
\mathbb{A}$ over $\mathbb{R}$, not necessarily associative or
finite-dimensional, equipped with a nonsingular quadratic form $Q:\mathbb{A}%
\rightarrow \mathbb{R}$ such that $Q(ab)=Q(a)Q(b)$ for all $a,b\in $ $%
\mathbb{A}$. The form $Q$ is given by the norm. For more information on
quadratic form, see [\ref{lewis2006}, pp. 44 and 53])

Briefly $\mathbb{H}$, the algebra of quaternions, has the following
properties:

\begin{itemize}
\item $\mathbb{H}$ is a $4$ dimensional non-commutative (Carley-Dickson)
algebra over the reals.

\item $\mathbb{H}$ is an associative algebra.

\item $\mathbb{H}$ is a division algebra, i.e. an algebra which is also a
division ring, i.e., each nonzero element of $\mathbb{H}$ is invertible.

\item $\mathbb{H}$ is a composition algebra.

\item $\mathbb{H}$ is a flexible algebra, i.e. $(pq)p=p(qp)$ for all $p,q\in 
\mathbb{H}$.

\item $\mathbb{H}$ is an alternative algebra, i.e. they have the property $%
p(pq)=(pp)q$ and $(qp)p=q(pp)$ for all $p,q\in \mathbb{H}$.
\end{itemize}

For the basics on the quaternions theory, we refer the work of Ward [\ref%
{ward1997}] and Lewis [\ref{lewis2006}].

We remark that

\begin{itemize}
\item $\mathbb{R}$, $\mathbb{C}$, $\mathbb{H}$ and $\mathbb{O}$ are the only
normed division algebras.

\item $\mathbb{R}$, $\mathbb{C}$, $\mathbb{H}$ and $\mathbb{O}$ are the only
alternative division algebras.
\end{itemize}

Last two properties shows what is so great about $\mathbb{R}$, $\mathbb{C}$, 
$\mathbb{H}$ and $\mathbb{O}$. For this two properties and their histories,
see [\ref{bib:baez2002}
].

\section{The Tetranacci and Tetranacci-Lucas Quaternions and their
Generating Functions, Binet's Formulas and Summations Formulas}

In this section, we define Tetranacci and Tetranacci-Lucas quaternions and
give generating functions and Binet formulas for them. First, we give some
information about quaternion sequences from the literature.

There are various types of quaternion sequences which have been studied by
many researchers. Horadam [\ref{bib:horadam1963aa}] introduced $n$th
Fibonacci and $n$th Lucas quaternions as%
\begin{equation*}
Q_{n}=F_{n}+F_{n+1}e_{1}+F_{n+2}e_{2}+F_{n+3}e_{3}=\sum_{s=0}^{3}F_{n+s}e_{s}
\end{equation*}%
and%
\begin{equation*}
R_{n}=L_{n}+L_{n+1}e_{1}+L_{n+2}e_{2}+L_{n+3}e_{3}=\sum_{s=0}^{3}L_{n+s}e_{s}
\end{equation*}%
respectively, where $F_{n}$ and $L_{n}$ are the $n$th Fibonacci and Lucas
numbers respectively. He also defined generalized Fibonacci quaternion as%
\begin{equation*}
P_{n}=H_{n}+H_{n+1}e_{1}+H_{n+2}e_{2}+H_{n+3}e_{3}=\sum_{s=0}^{3}H_{n+s}e_{s}
\end{equation*}%
where $H_{n}$ is the $n$th generalized Fibonacci number (which is now called
Horadam number) by the recursive relation $H_{1}=p,$ $H_{2}=p+q,$ $%
H_{n}=H_{n-1}+H_{n-2}$ ($p$ and $q$ are arbitrary integers). Halici [\ref%
{halici2012fibqua}] gave the generating functions and Binet formulas for the
Fibonacci and Lucas quaternions.

Cerda-Morales [\ref{bib:cerdamorale2017}] defined and studied the
generalized Tribonacci quaternion sequence that includes the previously
introduced Tribonacci, Padovan, Narayana and third order Jacobsthal
quaternion sequences. In [\ref{bib:cerdamorale2017}], the author defined
generalized Tribonacci quaternion as%
\begin{equation*}
Q_{v,n}=V_{n}+V_{n+1}e_{1}+V_{n+2}e_{2}+V_{n+3}e_{3}=%
\sum_{s=0}^{3}V_{n+s}e_{s}
\end{equation*}%
where $V_{n}$ is the $n$th generalized Tribonacci number defined by the
third-order recurrance relations%
\begin{equation*}
V_{n}=rV_{n-1}+sV_{n-2}+tV_{n-3},\text{ \ \ \ \ }
\end{equation*}%
here $V_{0}=a,V_{1}=b,V_{2}=c$ are arbitrary integers and $r,s,t$ are real
numbers.

Many other generalizations of Fibonacci quaternions have been given, see for
example Catarino [\ref{catarinopell2016}], Halici and Karata\c{s} [\ref%
{bib:halici2017}], and Polatl\i\ [\ref{bib:polatli2016}], Szynal-Liana and
Wloch [\ref{szynallianapell2016}] and Tasci [\ref{tascijacobquater2017}] for
second order quaternion sequences and Akkus and K\i z\i laslan [\ref%
{akkus2017}], Szynal-Liana and Wloch [\ref{szynal2017}], Tasci [\ref%
{tasciquater2018}], Cerda-Morales [\ref{cerdamoralejacobs2017}] for third
order quaternion sequences.

We now define Tetranacci and Tetranacci-Lucas quaternions over the
quaternion algebra $\mathbb{H}$. The $n$th Tetranacci quaternion is%
\begin{equation}
\widehat{M}_{n}=M_{n}+iM_{n+1}+jM_{n+2}+kM_{n+3}
\label{equation:bnuyvbcxdszaerm}
\end{equation}%
and the $n$th Tetranacci-Lucas quaternion is%
\begin{equation}
\widehat{R}_{n}=R_{n}+iR_{n+1}+jR_{n+2}+kR_{n+3}.
\end{equation}%
It can be easily shown that 
\begin{equation}
\widehat{M}_{n}=\widehat{M}_{n-1}+\widehat{M}_{n-2}+\widehat{M}_{n-3}+%
\widehat{M}_{n-4}  \label{equation:abgsvbnuytsaerds}
\end{equation}%
and%
\begin{equation}
\widehat{R}_{n}=\widehat{R}_{n-1}+\widehat{R}_{n-2}+\widehat{R}_{n-3}+%
\widehat{R}_{n-4}.  \label{equation:cvnbhustarfcxdsa}
\end{equation}%
Note that 
\begin{equation*}
\widehat{M}_{-n}=-\widehat{M}_{-(n-1)}-\widehat{M}_{-(n-2)}-\widehat{M}%
_{-(n-3)}+\widehat{M}_{-(n-4)}
\end{equation*}%
and%
\begin{equation*}
\widehat{R}_{-n}=-\widehat{R}_{-(n-1)}-\widehat{R}_{-(n-2)}-\widehat{R}%
_{-(n-3)}+\widehat{R}_{-(n-4)}.
\end{equation*}

The conjugate of $\widehat{M}_{n}$ and $\widehat{R}_{n}$ are defined by%
\begin{equation*}
\overline{\widehat{M}_{n}}=M_{n}-iM_{n+1}-jM_{n+2}-kM_{n+3}
\end{equation*}%
and%
\begin{equation*}
\overline{\widehat{R}_{n}}=R_{n}-iR_{n+1}-jR_{n+2}-kR_{n+3}
\end{equation*}%
respectively.

Now, we will state Binet's formula for the Tetranacci and Tetranacci-Lucas
quaternions and in the rest of the paper we fix the following notations.%
\begin{eqnarray*}
\widehat{\alpha } &=&1+i\alpha +j\alpha ^{2}+k\alpha ^{3}, \\
\widehat{\beta } &=&1+i\beta +j\beta ^{2}+k\beta ^{3}, \\
\widehat{\gamma } &=&1+i\gamma +j\gamma ^{2}+k\gamma ^{3}. \\
\widehat{\delta } &=&1+i\delta +j\delta ^{2}+k\delta ^{3}.
\end{eqnarray*}

\begin{theorem}
\label{theorem:yunhbvgcfosea}(Binet's Formulas) For any integer $n,$ the $n$%
th Tetranacci quaternion is%
\begin{eqnarray}
\widehat{M}_{n} &=&\frac{\widehat{\alpha }\alpha ^{n+2}}{(\alpha -\beta
)(\alpha -\gamma )(\alpha -\delta )}+\frac{\widehat{\beta }\beta ^{n+2}}{%
(\beta -\alpha )(\beta -\gamma )(\beta -\delta )}
\label{equation:fgvbvcxsdaewqzxd} \\
&&+\frac{\widehat{\gamma }\gamma ^{n+2}}{(\gamma -\alpha )(\gamma -\beta
)(\gamma -\delta )}+\frac{\widehat{\delta }\delta ^{n+2}}{(\delta -\alpha
)(\delta -\beta )(\delta -\gamma )}  \notag \\
&=&\frac{\alpha -1}{5\alpha -8}\widehat{\alpha }\alpha ^{n-1}+\frac{\beta -1%
}{5\beta -8}\widehat{\beta }\beta ^{n-1}+\frac{\gamma -1}{5\gamma -8}%
\widehat{\gamma }\gamma ^{n-1}+\frac{\delta -1}{5\delta -8}\widehat{\delta }%
\delta ^{n-1}  \label{equati:gbvgcdxscazuyt}
\end{eqnarray}%
and the $n$th Tetranacci-Lucas quaternion is 
\begin{equation}
\widehat{R}_{n}=\widehat{\alpha }\alpha ^{n}+\widehat{\beta }\beta ^{n}+%
\widehat{\gamma }\gamma ^{n}+\widehat{\delta }\delta ^{n}.
\label{equation:xcdfvesaztfvop}
\end{equation}
\end{theorem}

\textit{Proof.} Using Binet's formula of the Tetranacci-Lucas numbers, we
have 
\begin{eqnarray*}
\widehat{R}_{n} &=&R_{n}+iR_{n+1}+jR_{n+2}+kR_{n+3} \\
&=&(\alpha ^{n}+\beta ^{n}+\gamma ^{n}+\delta ^{n})+i(\alpha ^{n+1}+\beta
^{n+1}+\gamma ^{n+1}+\delta ^{n+1}) \\
&&+j(\alpha ^{n+2}+\beta ^{n+2}+\gamma ^{n+2}+\delta ^{n+2})+k(\alpha
^{n+3}+\beta ^{n+3}+\gamma ^{n+3}+\delta ^{n+3}) \\
&=&\widehat{\alpha }\alpha ^{n}+\widehat{\beta }\beta ^{n}+\widehat{\gamma }%
\gamma ^{n}+\widehat{\delta }\delta ^{n}.
\end{eqnarray*}

Note that using Binet's formula (\ref{equation:srtynsbfdtsra}) of the
Tetranacci numbers we have%
\begin{eqnarray*}
\widehat{M}_{n} &=&M_{n}+iM_{n+1}+jM_{n+2}+kM_{n+3} \\
&=&(\frac{\alpha -1}{5\alpha -8}\alpha ^{n-1}+\frac{\beta -1}{5\beta -8}%
\beta ^{n-1}+\frac{\gamma -1}{5\gamma -8}\gamma ^{n-1}+\frac{\delta -1}{%
5\delta -8}\delta ^{n-1}) \\
&&+i(\frac{\alpha -1}{5\alpha -8}\alpha ^{n}+\frac{\beta -1}{5\beta -8}\beta
^{n}+\frac{\gamma -1}{5\gamma -8}\gamma ^{n}+\frac{\delta -1}{5\delta -8}%
\delta ^{n}) \\
&&+j(\frac{\alpha -1}{5\alpha -8}\alpha ^{n+1}+\frac{\beta -1}{5\beta -8}%
\beta ^{n+1}+\frac{\gamma -1}{5\gamma -8}\gamma ^{n+1}+\frac{\delta -1}{%
5\delta -8}\delta ^{n+1}) \\
&&+k(\frac{\alpha -1}{5\alpha -8}\alpha ^{n+2}+\frac{\beta -1}{5\beta -8}%
\beta ^{n+2}+\frac{\gamma -1}{5\gamma -8}\gamma ^{n+2}+\frac{\delta -1}{%
5\delta -8}\delta ^{n+2}) \\
&=&\frac{\alpha -1}{5\alpha -8}\widehat{\alpha }\alpha ^{n-1}+\frac{\beta -1%
}{5\beta -8}\widehat{\beta }\beta ^{n-1}+\frac{\gamma -1}{5\gamma -8}%
\widehat{\gamma }\gamma ^{n-1}+\frac{\delta -1}{5\delta -8}\widehat{\delta }%
\delta ^{n-1}.
\end{eqnarray*}%
This proves (\ref{equati:gbvgcdxscazuyt}). Similarly, we can obtain (\ref%
{equation:fgvbvcxsdaewqzxd}). 
\endproof%

\begin{remark}
According to Theorem \ref{theorem:fvgxdfsxczsaer}, Binet's Formulas of the
Tetranacci and Tetranacci-Lucas quaternions are true for all integers $n.$
\end{remark}

Next, we present generating functions.

\begin{theorem}
The generating functions for the Tetranacci and Tetranacci-Lucas quaternions
are%
\begin{equation}
\sum_{n=0}^{\infty }\widehat{M}_{n}x^{n}=\frac{%
(i+j+2k)+(1+j+2k)x+(j+2k)x^{2}+(j+k)x^{3}}{1-x-x^{2}-x^{3}-x^{4}}
\label{equation:hnbvyuedsxzapo}
\end{equation}%
and%
\begin{equation}
\sum_{n=0}^{\infty }\widehat{R}_{n}x^{n}=\frac{%
(4+i+3j+7k)+(-3+2i+4j+8k)x+(-2+3i+5j+4k)x^{2}+(-1+4i+j+3k)x^{3}}{%
1-x-x^{2}-x^{3}-x^{4}}  \label{equat:bnopscxzartfnm}
\end{equation}%
respectively.
\end{theorem}

\textit{Proof.} Let%
\begin{equation*}
g(x)=\sum_{n=0}^{\infty }\widehat{M}_{n}x^{n}
\end{equation*}%
be generating function of\ the Tetranacci quaternions. Then, using the
definition of the Tetranacci quaternions, and substracting $xg(x),$ $%
x^{2}g(x),$ $x^{3}g(x)$ and $x^{4}g(x)$ from $g(x),$ we obtain (note the
shift in the index $n$ in the third line)%
\begin{eqnarray*}
&&(1-x-x^{2}-x^{3}-x^{4})g(x) \\
&=&\sum_{n=0}^{\infty }\widehat{M}_{n}x^{n}-x\sum_{n=0}^{\infty }\widehat{M}%
_{n}x^{n}-x^{2}\sum_{n=0}^{\infty }\widehat{M}_{n}x^{n}-x^{3}\sum_{n=0}^{%
\infty }\widehat{M}_{n}x^{n}-x^{4}\sum_{n=0}^{\infty }\widehat{M}_{n}x^{n} \\
&=&\sum_{n=0}^{\infty }\widehat{M}_{n}x^{n}-\sum_{n=0}^{\infty }\widehat{M}%
_{n}x^{n+1}-\sum_{n=0}^{\infty }\widehat{M}_{n}x^{n+2}-\sum_{n=0}^{\infty }%
\widehat{M}_{n}x^{n+3}-\sum_{n=0}^{\infty }\widehat{M}_{n}x^{n+4} \\
&=&\sum_{n=0}^{\infty }\widehat{M}_{n}x^{n}-\sum_{n=1}^{\infty }\widehat{M}%
_{n-1}x^{n}-\sum_{n=2}^{\infty }\widehat{M}_{n-2}x^{n}-\sum_{n=3}^{\infty }%
\widehat{M}_{n-3}x^{n}-\sum_{n=4}^{\infty }\widehat{M}_{n-4}x^{n} \\
&=&(\widehat{M}_{0}+\widehat{M}_{1}x+\widehat{M}_{2}x^{2}+\widehat{M}%
_{3}x^{3})-(\widehat{M}_{0}x+\widehat{M}_{1}x^{2}+\widehat{M}_{2}x^{3})-(%
\widehat{M}_{0}x^{2}+\widehat{M}_{1}x^{3})-\widehat{M}_{0}x^{3} \\
&&+\sum_{n=4}^{\infty }(\widehat{M}_{n}-\widehat{M}_{n-1}-\widehat{M}_{n-2}-%
\widehat{M}_{n-3}-\widehat{M}_{n-4})x^{n} \\
&=&\widehat{M}_{0}+(\widehat{M}_{1}-\widehat{M}_{0})x+(\widehat{M}_{2}-%
\widehat{M}_{1}-\widehat{M}_{0})x^{2}+(\widehat{M}_{3}-\widehat{M}_{2}-%
\widehat{M}_{1}-\widehat{M}_{0})x^{3}.
\end{eqnarray*}%
Note that we used the recurrence relation $\widehat{M}_{n}=\widehat{M}_{n-1}+%
\widehat{M}_{n-2}+\widehat{M}_{n-3}+\widehat{M}_{n-4}$. Rearranging above
equation, we get%
\begin{equation*}
g(x)=\frac{\widehat{M}_{0}+(\widehat{M}_{1}-\widehat{M}_{0})x+(\widehat{M}%
_{2}-\widehat{M}_{1}-\widehat{M}_{0})x^{2}+(\widehat{M}_{3}-\widehat{M}_{2}-%
\widehat{M}_{1}-\widehat{M}_{0})x^{3}}{1-x-x^{2}-x^{3}-x^{4}}.
\end{equation*}%
or 
\begin{equation*}
g(x)=\frac{\widehat{M}_{0}+(\widehat{M}_{1}-\widehat{M}_{0})x+(\widehat{M}%
_{2}-\widehat{M}_{1}-\widehat{M}_{0})x^{2}+\widehat{M}_{-1}x^{3}}{%
1-x-x^{2}-x^{3}-x^{4}}.
\end{equation*}%
since $\widehat{M}_{3}=\widehat{M}_{2}+\widehat{M}_{1}+\widehat{M}_{0}+%
\widehat{M}_{-1}$. Now using%
\begin{eqnarray*}
\widehat{M}_{-1} &=&j+k, \\
\widehat{M}_{0} &=&i+j+2k, \\
\widehat{M}_{1} &=&1+i+2j+4k, \\
\widehat{M}_{2} &=&1+2i+4j+8k, \\
\widehat{M}_{3} &=&2+4i+8j+15k,
\end{eqnarray*}%
we obtain%
\begin{equation*}
g(x)=\frac{(i+j+2k)+(1+j+2k)x+(j+2k)x^{2}+(j+k)x^{3}}{1-x-x^{2}-x^{3}-x^{4}}.
\end{equation*}%
Similarly, we can obtain (\ref{equat:bnopscxzartfnm}). 
\endproof%

In the following theorem, we present another forms of Binet's formulas for
the Tetranacci and Tetranacci-Lucas quaternions using generating functions.

\begin{theorem}
\label{theorem:vbghcxdszrtupomn}For any integer $n,$ the $n$th Tetranacci
quaternion is%
\begin{eqnarray*}
\widehat{M}_{n} &=&\frac{\widehat{M}_{-1}+\alpha (\widehat{M}_{2}-\widehat{M}%
_{1}-\widehat{M}_{0})+\alpha ^{2}\allowbreak (\widehat{M}_{1}-\widehat{M}%
_{0})+\alpha ^{3}\widehat{M}_{0}}{\left( \alpha -\beta \right) \left( \alpha
-\gamma \right) \left( \alpha -\delta \right) }\alpha ^{n} \\
&&+\frac{\widehat{M}_{-1}+\beta (\widehat{M}_{2}-\widehat{M}_{1}-\widehat{M}%
_{0})+\beta ^{2}\allowbreak (\widehat{M}_{1}-\widehat{M}_{0})+\beta ^{3}%
\widehat{M}_{0}}{\left( \beta -\gamma \right) \left( \beta -\alpha \right)
\left( \beta -\delta \right) }\beta ^{n} \\
&&+\frac{\widehat{M}_{-1}+\gamma (\widehat{M}_{2}-\widehat{M}_{1}-\widehat{M}%
_{0})+\gamma ^{2}\allowbreak (\widehat{M}_{1}-\widehat{M}_{0})+\gamma ^{3}%
\widehat{M}_{0}}{(\gamma -\alpha )(\gamma -\beta )(\gamma -\delta )}\gamma
^{n} \\
&&+\frac{\widehat{M}_{-1}+\delta (\widehat{M}_{2}-\widehat{M}_{1}-\widehat{M}%
_{0})+\delta ^{2}\allowbreak (\widehat{M}_{1}-\widehat{M}_{0})+\delta ^{3}%
\widehat{M}_{0}}{(\delta -\alpha )(\delta -\beta )(\delta -\gamma )}\delta
^{n}
\end{eqnarray*}%
and the $n$th Tetranacci-Lucas quaternion is%
\begin{eqnarray*}
\widehat{R}_{n} &=&\frac{\widehat{R}_{-1}+\alpha (\widehat{R}_{2}-\widehat{R}%
_{1}-\widehat{R}_{0})+\alpha ^{2}\allowbreak (\widehat{R}_{1}-\widehat{R}%
_{0})+\alpha ^{3}\widehat{R}_{0}}{\left( \alpha -\beta \right) \left( \alpha
-\gamma \right) \left( \alpha -\delta \right) }\alpha ^{n}+\frac{\widehat{R}%
_{-1}+\beta (\widehat{R}_{2}-\widehat{R}_{1}-\widehat{R}_{0})+\beta
^{2}\allowbreak (\widehat{R}_{1}-\widehat{R}_{0})+\beta ^{3}\widehat{R}_{0}}{%
\left( \beta -\gamma \right) \left( \beta -\alpha \right) \left( \beta
-\delta \right) }\beta ^{n} \\
&&+\frac{\widehat{R}_{-1}+\gamma (\widehat{R}_{2}-\widehat{R}_{1}-\widehat{R}%
_{0})+\gamma ^{2}\allowbreak (\widehat{R}_{1}-\widehat{R}_{0})+\gamma ^{3}%
\widehat{R}_{0}}{(\gamma -\alpha )(\gamma -\beta )(\gamma -\delta )}\gamma
^{n}+\frac{\widehat{R}_{-1}+\delta (\widehat{R}_{2}-\widehat{R}_{1}-\widehat{%
R}_{0})+\delta ^{2}\allowbreak (\widehat{R}_{1}-\widehat{R}_{0})+\delta ^{3}%
\widehat{R}_{0}}{(\delta -\alpha )(\delta -\beta )(\delta -\gamma )}\delta
^{n}
\end{eqnarray*}
\end{theorem}

\textit{Proof.} We can use generating functions. Since the roots of the
equation $1-x-x^{2}-x^{3}-x^{4}=0$ are $\frac{1}{\alpha },\frac{1}{\beta },%
\frac{1}{\gamma },\frac{1}{\delta }$ and%
\begin{equation*}
1-x-x^{2}-x^{3}-x^{4}=(1-\alpha x)(1-\beta x)(1-\gamma x)(1-\delta x),
\end{equation*}%
we can write the generating function of $\widehat{M}_{n}$ as%
\begin{eqnarray*}
g(x) &=&\frac{\widehat{M}_{0}+(\widehat{M}_{1}-\widehat{M}_{0})x+(\widehat{M}%
_{2}-\widehat{M}_{1}-\widehat{M}_{0})x^{2}+\widehat{M}_{-1}x^{3}}{%
1-x-x^{2}-x^{3}-x^{4}} \\
&=&\frac{\widehat{M}_{0}+(\widehat{M}_{1}-\widehat{M}_{0})x+(\widehat{M}_{2}-%
\widehat{M}_{1}-\widehat{M}_{0})x^{2}+\widehat{M}_{-1}x^{3}}{(1-\alpha
x)(1-\beta x)(1-\gamma x)(1-\delta x)} \\
&=&\frac{A}{(1-\alpha x)}+\frac{B}{(1-\beta x)}+\frac{C}{(1-\gamma x)}+\frac{%
D}{(1-\delta x)}
\end{eqnarray*}%
We need to find $A,B,C$ and $D$, so the following system of equations should
be solved:%
\begin{eqnarray*}
A+B+C+D &=&\widehat{M}_{0} \\
A(-\beta -\gamma -\delta )+B(-\alpha -\gamma -\delta )+C(-\alpha -\beta
-\delta )+D(-\alpha -\beta -\gamma ) &=&\widehat{M}_{1}-\widehat{M}_{0} \\
A(\beta \gamma +\beta \delta +\gamma \allowbreak \delta )+B(\alpha \gamma
+\alpha \delta +\gamma \allowbreak \delta )+C(\alpha \beta +\alpha \delta
+\beta \allowbreak \delta )+D(\alpha \beta +\alpha \gamma +\beta \gamma
\allowbreak ) &=&\widehat{M}_{2}-\widehat{M}_{1}-\widehat{M}_{0} \\
-A\beta \gamma \delta -B\alpha \gamma \delta -C\alpha \beta \delta -\alpha
\beta \gamma D &=&\widehat{M}_{-1}.
\end{eqnarray*}%
Then, we find that%
\begin{eqnarray*}
A &=&\frac{\widehat{M}_{-1}+\alpha (\widehat{M}_{2}-\widehat{M}_{1}-\widehat{%
M}_{0})+\alpha ^{2}\allowbreak (\widehat{M}_{1}-\widehat{M}_{0})+\alpha ^{3}%
\widehat{M}_{0}}{\left( \alpha -\beta \right) \left( \alpha -\gamma \right)
\left( \alpha -\delta \right) } \\
B &=&\frac{\widehat{M}_{-1}+\beta (\widehat{M}_{2}-\widehat{M}_{1}-\widehat{M%
}_{0})+\beta ^{2}\allowbreak (\widehat{M}_{1}-\widehat{M}_{0})+\beta ^{3}%
\widehat{M}_{0}}{\left( \beta -\gamma \right) \left( \beta -\alpha \right)
\left( \beta -\delta \right) } \\
C &=&\frac{\widehat{M}_{-1}+\gamma (\widehat{M}_{2}-\widehat{M}_{1}-\widehat{%
M}_{0})+\gamma ^{2}\allowbreak (\widehat{M}_{1}-\widehat{M}_{0})+\gamma ^{3}%
\widehat{M}_{0}}{(\gamma -\alpha )(\gamma -\beta )(\gamma -\delta )} \\
D &=&\frac{\widehat{M}_{-1}+\delta (\widehat{M}_{2}-\widehat{M}_{1}-\widehat{%
M}_{0})+\delta ^{2}\allowbreak (\widehat{M}_{1}-\widehat{M}_{0})+\delta ^{3}%
\widehat{M}_{0}}{(\delta -\alpha )(\delta -\beta )(\delta -\gamma )}
\end{eqnarray*}%
and 
\begin{eqnarray*}
g(x) &=&\frac{\widehat{M}_{-1}+\alpha (\widehat{M}_{2}-\widehat{M}_{1}-%
\widehat{M}_{0})+\alpha ^{2}\allowbreak (\widehat{M}_{1}-\widehat{M}%
_{0})+\alpha ^{3}\widehat{M}_{0}}{\left( \alpha -\beta \right) \left( \alpha
-\gamma \right) \left( \alpha -\delta \right) }\sum_{n=0}^{\infty }\alpha
^{n}x^{n} \\
&&+\frac{\widehat{M}_{-1}+\beta (\widehat{M}_{2}-\widehat{M}_{1}-\widehat{M}%
_{0})+\beta ^{2}\allowbreak (\widehat{M}_{1}-\widehat{M}_{0})+\beta ^{3}%
\widehat{M}_{0}}{\left( \beta -\gamma \right) \left( \beta -\alpha \right)
\left( \beta -\delta \right) }\sum_{n=0}^{\infty }\beta ^{n}x^{n} \\
&&+\frac{\widehat{M}_{-1}+\gamma (\widehat{M}_{2}-\widehat{M}_{1}-\widehat{M}%
_{0})+\gamma ^{2}\allowbreak (\widehat{M}_{1}-\widehat{M}_{0})+\gamma ^{3}%
\widehat{M}_{0}}{(\gamma -\alpha )(\gamma -\beta )(\gamma -\delta )}%
\sum_{n=0}^{\infty }\gamma ^{n}x^{n} \\
&&+\frac{\widehat{M}_{-1}+\delta (\widehat{M}_{2}-\widehat{M}_{1}-\widehat{M}%
_{0})+\delta ^{2}\allowbreak (\widehat{M}_{1}-\widehat{M}_{0})+\delta ^{3}%
\widehat{M}_{0}}{(\delta -\alpha )(\delta -\beta )(\delta -\gamma )}%
\sum_{n=0}^{\infty }\delta ^{n}x^{n} \\
&=&\sum_{n=0}^{\infty }\left( 
\begin{array}{c}
\frac{\widehat{M}_{-1}+\alpha (\widehat{M}_{2}-\widehat{M}_{1}-\widehat{M}%
_{0})+\alpha ^{2}\allowbreak (\widehat{M}_{1}-\widehat{M}_{0})+\alpha ^{3}%
\widehat{M}_{0}}{\left( \alpha -\beta \right) \left( \alpha -\gamma \right)
\left( \alpha -\delta \right) }\alpha ^{n}+\frac{\widehat{M}_{-1}+\beta (%
\widehat{M}_{2}-\widehat{M}_{1}-\widehat{M}_{0})+\beta ^{2}\allowbreak (%
\widehat{M}_{1}-\widehat{M}_{0})+\beta ^{3}\widehat{M}_{0}}{\left( \beta
-\gamma \right) \left( \beta -\alpha \right) \left( \beta -\delta \right) }%
\beta ^{n} \\ 
+\frac{\widehat{M}_{-1}+\gamma (\widehat{M}_{2}-\widehat{M}_{1}-\widehat{M}%
_{0})+\gamma ^{2}\allowbreak (\widehat{M}_{1}-\widehat{M}_{0})+\gamma ^{3}%
\widehat{M}_{0}}{(\gamma -\alpha )(\gamma -\beta )(\gamma -\delta )}\gamma
^{n}+\frac{\widehat{M}_{-1}+\delta (\widehat{M}_{2}-\widehat{M}_{1}-\widehat{%
M}_{0})+\delta ^{2}\allowbreak (\widehat{M}_{1}-\widehat{M}_{0})+\delta ^{3}%
\widehat{M}_{0}}{(\delta -\alpha )(\delta -\beta )(\delta -\gamma )}\delta
^{n}%
\end{array}%
\right) x^{n}.
\end{eqnarray*}%
Thus, from this, we obtain Binet's formula of Tetranacci quaternion.
Similarly, we can obtain Binet's formula of the Tetranacci-Lucas quaternion. 
\endproof%

If we compare Theorem \ref{theorem:yunhbvgcfosea} and Theorem \ref%
{theorem:vbghcxdszrtupomn} and use the definition of $\widehat{M}_{n},$ $%
\widehat{R}_{n},$ we have the following Remark showing relations between $%
\widehat{M}_{-1},\widehat{M}_{0},\widehat{M}_{1},\widehat{M}_{2};\widehat{R}%
_{-1},\widehat{R}_{0},\widehat{R}_{1},\widehat{R}_{2}$ and $\widehat{\alpha }%
,\widehat{\beta },\widehat{\gamma },\widehat{\delta }.$

\begin{remark}
We have the following identities:

\begin{description}
\item[(a)] 
\begin{eqnarray*}
\frac{\widehat{M}_{-1}+\alpha (\widehat{M}_{2}-\widehat{M}_{1}-\widehat{M}%
_{0})+\alpha ^{2}\allowbreak (\widehat{M}_{1}-\widehat{M}_{0})+\alpha ^{3}%
\widehat{M}_{0}}{\alpha ^{2}} &=&\widehat{\alpha } \\
\frac{\widehat{M}_{-1}+\beta (\widehat{M}_{2}-\widehat{M}_{1}-\widehat{M}%
_{0})+\beta ^{2}\allowbreak (\widehat{M}_{1}-\widehat{M}_{0})+\beta ^{3}%
\widehat{M}_{0}}{\beta ^{2}} &=&\widehat{\beta } \\
\frac{\widehat{M}_{-1}+\gamma (\widehat{M}_{2}-\widehat{M}_{1}-\widehat{M}%
_{0})+\gamma ^{2}\allowbreak (\widehat{M}_{1}-\widehat{M}_{0})+\gamma ^{3}%
\widehat{M}_{0}}{\gamma ^{2}} &=&\widehat{\gamma } \\
\frac{\widehat{M}_{-1}+\delta (\widehat{M}_{2}-\widehat{M}_{1}-\widehat{M}%
_{0})+\delta ^{2}\allowbreak (\widehat{M}_{1}-\widehat{M}_{0})+\delta ^{3}%
\widehat{M}_{0}}{\delta ^{2}} &=&\widehat{\delta }
\end{eqnarray*}

\item[(b)] 
\begin{eqnarray*}
\frac{\widehat{R}_{-1}+\alpha (\widehat{R}_{2}-\widehat{R}_{1}-\widehat{R}%
_{0})+\alpha ^{2}\allowbreak (\widehat{R}_{1}-\widehat{R}_{0})+\alpha ^{3}%
\widehat{R}_{0}}{\left( \alpha -\beta \right) \left( \alpha -\gamma \right)
\left( \alpha -\delta \right) } &=&\widehat{\alpha } \\
+\frac{\widehat{R}_{-1}+\beta (\widehat{R}_{2}-\widehat{R}_{1}-\widehat{R}%
_{0})+\beta ^{2}\allowbreak (\widehat{R}_{1}-\widehat{R}_{0})+\beta ^{3}%
\widehat{R}_{0}}{\left( \beta -\gamma \right) \left( \beta -\alpha \right)
\left( \beta -\delta \right) } &=&\widehat{\beta } \\
\frac{\widehat{R}_{-1}+\gamma (\widehat{R}_{2}-\widehat{R}_{1}-\widehat{R}%
_{0})+\gamma ^{2}\allowbreak (\widehat{R}_{1}-\widehat{R}_{0})+\gamma ^{3}%
\widehat{R}_{0}}{(\gamma -\alpha )(\gamma -\beta )(\gamma -\delta )} &=&%
\widehat{\gamma } \\
\frac{\widehat{R}_{-1}+\delta (\widehat{R}_{2}-\widehat{R}_{1}-\widehat{R}%
_{0})+\delta ^{2}\allowbreak (\widehat{R}_{1}-\widehat{R}_{0})+\delta ^{3}%
\widehat{R}_{0}}{(\delta -\alpha )(\delta -\beta )(\delta -\gamma )} &=&%
\widehat{\delta }
\end{eqnarray*}
\end{description}
\end{remark}

Now, we present the formulas which give the summation of the first $n$
Tetranacci and Tetranacci-Lucas numbers.

\begin{lemma}
For every integer $n\geq 0,$ we have 
\begin{equation}
\sum\limits_{p=0}^{n}M_{p}=\frac{1}{3}(M_{n+2}+2M_{n}+M_{n-1}-1)
\label{equation:easzpuyhbnvcdf}
\end{equation}%
and%
\begin{equation}
\sum\limits_{p=0}^{n}R_{p}=\frac{1}{3}(R_{n+2}+2R_{n}+R_{n-1}+2).
\label{equation:prtsaeoaxcdf}
\end{equation}
\end{lemma}

\textit{Proof.} (\ref{equation:easzpuyhbnvcdf}) and (\ref%
{equation:prtsaeoaxcdf}) are given in Soykan [\ref%
{soykan2019gaussiangentetran}, Corollaries 2.7 and 2.8]. 
\endproof%

Note that (\ref{equation:easzpuyhbnvcdf}) and (\ref{equation:prtsaeoaxcdf})
can be easily proved by mathematical induction as well.

Next, we present the formulas which give the summation of the first $n$
Tetranacci and Tetranacci-Lucas quaternions.

\begin{theorem}
The summation formula for Tetranacci and Tetranacci-Lucas quaternions are%
\begin{equation}
\sum\limits_{p=0}^{n}\widehat{M}_{p}=\frac{1}{3}(\widehat{M}_{n+2}+2\widehat{%
M}_{n}+\widehat{M}_{n-1}-(1+i+4j+7k))  \label{equation:yaomazxtyhnbvgfc}
\end{equation}%
and%
\begin{equation}
\sum\limits_{p=0}^{n}\widehat{R}_{p}=\frac{1}{3}(\widehat{R}_{n+2}+2\widehat{%
R}_{n}+\widehat{R}_{n-1}+(2-10i-13j-22k)).  \label{equation:mguouresxcdz}
\end{equation}
\end{theorem}

\textit{Proof.} Using (\ref{equation:bnuyvbcxdszaerm}) and (\ref%
{equation:easzpuyhbnvcdf}), we obtain%
\begin{eqnarray*}
\sum\limits_{p=0}^{n}\widehat{M}_{i}
&=&\sum\limits_{p=0}^{n}M_{p}+i\sum\limits_{p=0}^{n}M_{p+1}+j\sum%
\limits_{p=0}^{n}M_{p+2}+k\sum\limits_{p=0}^{n}M_{p+3} \\
&=&(M_{0}+...+M_{n})+i(M_{1}+...+M_{n+1}) \\
&&+j(M_{2}+...+M_{n+2})+k(M_{3}+...+M_{n+3}).
\end{eqnarray*}%
and so%
\begin{eqnarray*}
3\sum\limits_{p=0}^{n}\widehat{M}_{p} &=&(M_{n+2}+2M_{n}+M_{n-1}-1) \\
&&+i(M_{n+3}+2M_{n+1}+M_{n}-1-3M_{0}) \\
&&+j(M_{n+4}+2M_{n+2}+M_{n+1}-1-3(M_{0}+M_{1})) \\
&&+k(M_{n+5}+2M_{n+3}+M_{n+2}-1-3(M_{0}+M_{1}+M_{2})) \\
&=&\widehat{M}_{n+2}+2\widehat{M}_{n}+\widehat{M}_{n-1}+c
\end{eqnarray*}

where 
\begin{eqnarray*}
c &=&-1+i(-1-3M_{0})+j(-1-3(M_{0}+M_{1}))+k(-1-3(M_{0}+M_{1}+M_{2})) \\
&=&-1-i-4j-7k.
\end{eqnarray*}%
Hence 
\begin{equation*}
\sum\limits_{p=0}^{n}\widehat{M}_{p}=\frac{1}{3}(\widehat{M}_{n+2}+2\widehat{%
M}_{n}+\widehat{M}_{n-1}-(1+i+4j+7k)).
\end{equation*}

This proves (\ref{equation:yaomazxtyhnbvgfc}). Similarly, we can obtain (\ref%
{equation:mguouresxcdz}). 
\endproof%

Note that above Theorem can be proved by induction as well.

\begin{theorem}
For $n\geq 0,$ we have the following formulas:

\begin{description}
\item[(a)] $\sum\limits_{p=0}^{n}\widehat{M}_{2p+1}=\frac{1}{3}(2\widehat{M}%
_{2n+2}+\widehat{M}_{2n}-\widehat{M}_{2n-1}+(1-2i-2j-5k))$

\item[(b)] $\sum\limits_{p=0}^{n}\widehat{M}_{2p}=\frac{1}{3}(2\widehat{M}%
_{2n+1}+\widehat{M}_{2n-1}-\widehat{M}_{2n-2}-(2-i+2j+2k)).$
\end{description}
\end{theorem}

\textit{Proof. }The proof follows from the following identities:%
\begin{equation}
\sum\limits_{p=0}^{n}M_{2p+1}=\frac{1}{3}(2M_{2n+2}+M_{2n}-M_{2n-1}+1)
\label{equati:hgfdserat}
\end{equation}%
and 
\begin{equation}
\sum\limits_{p=0}^{n}M_{2p}=\frac{1}{3}(2M_{2n+1}+M_{2n-1}-M_{2n-2}-2).
\label{equat:pcxzdsogdsap}
\end{equation}%
(\ref{equati:hgfdserat}) and (\ref{equat:pcxzdsogdsap}) are given in Soykan [%
\ref{soykan2019gaussiangentetran}, Corollary 2.7].\ 
\endproof%

Note that (\ref{equati:hgfdserat}) and (\ref{equat:pcxzdsogdsap}) can be
easily proved by mathematical induction as well. Of course, the above
theorem itself can be proved by induction.\ 

\begin{theorem}
For $n\geq 0,$ we have the following formulas:

\begin{description}
\item[(a)] $\sum\limits_{p=0}^{n}\widehat{R}_{2p+1}=\frac{1}{3}(2\widehat{R}%
_{2n+2}+\widehat{R}_{2n}-\widehat{R}_{2n-1}-(8+2i+11j+11k))$

\item[(b)] $\sum\limits_{p=0}^{n}\widehat{R}_{2p}=\frac{1}{3}(2\widehat{R}%
_{2n+1}+\widehat{R}_{2n-1}-\widehat{R}_{2n-2}+(10-8i-2j-11k)).$
\end{description}
\end{theorem}

\textit{Proof. }The proof follows from the following identities:%
\begin{equation}
\sum\limits_{p=0}^{n}R_{2p+1}=\frac{1}{3}(2R_{2n+2}+R_{2n}-R_{2n-1}-8)
\label{equati:kshpsmnbx}
\end{equation}%
and 
\begin{equation}
\sum\limits_{p=0}^{n}R_{2p}=\frac{1}{3}(2R_{2n+1}+R_{2n-1}-R_{2n-2}+10).
\label{equation:xcxdzsaopuytrs}
\end{equation}%
(\ref{equati:kshpsmnbx}) and (\ref{equation:xcxdzsaopuytrs}) are given in
Soykan [\ref{soykan2019gaussiangentetran}, Corollary 2.8].\ 
\endproof%

Note that (\ref{equati:kshpsmnbx}) and (\ref{equation:xcxdzsaopuytrs}) can
be easily proved by mathematical induction as well. Of course, the above
theorem itself can be proved by induction.\ 

\section{Matrices and Determinants related with Tetranacci and
Tetranacci-Lucas Quaternions}

Define the $5\times 5$ determinants $D_{n}$ and $E_{n},$ for all integers $n,
$ by%
\begin{equation*}
D_{n}=\left\vert 
\begin{array}{ccccc}
M_{n} & R_{n} & R_{n+1} & R_{n+2} & R_{n+3} \\ 
M_{2} & R_{2} & R_{3} & R_{4} & R_{5} \\ 
M_{1} & R_{1} & R_{2} & R_{3} & R_{4} \\ 
M_{0} & R_{0} & R_{1} & R_{2} & R_{3} \\ 
M_{-1} & R_{-1} & R_{0} & R_{1} & R_{2}%
\end{array}%
\right\vert ,\text{ }E_{n}=\left\vert 
\begin{array}{ccccc}
R_{n} & M_{n} & M_{n+1} & M_{n+2} & M_{n+3} \\ 
R_{2} & M_{2} & M_{3} & M_{4} & M_{5} \\ 
R_{1} & M_{1} & M_{2} & M_{3} & M_{4} \\ 
R_{0} & M_{0} & M_{1} & M_{2} & M_{3} \\ 
R_{-1} & M_{-1} & M_{0} & M_{1} & M_{2}%
\end{array}%
\right\vert .
\end{equation*}

\begin{theorem}
The following statements are true.

\begin{description}
\item[(a)] $D_{n}=0$ and $E_{n}=0$ for all integers $n.$

\item[(b)] $563\widehat{M}_{n}=86\widehat{R}_{n+3}-61\widehat{R}_{n+2}-71%
\widehat{R}_{n+1}-87\widehat{R}_{n}.$

\item[(c)] $\widehat{R}_{n}=6\widehat{M}_{n+1}-\widehat{M}_{n}-\widehat{M}%
_{n+3}.$
\end{description}
\end{theorem}

\textit{Proof.} (a) is a special case of a result in [\ref{bib:melham1995}].
Expanding $D_{n}$ along the top row gives $%
563M_{n}=86R_{n+3}-61R_{n+2}-71R_{n+1}-87R_{n}\ $and now (b) follows.
Expanding $E_{n}$ along the top row gives $R_{n}=6M_{n+1}-M_{n}-M_{n+3}\ $%
and now (c) follows. 
\endproof%

Consider the sequence $\{U_{n}\}$ which is defined by the fourth-order
recurrence relation%
\begin{equation*}
U_{n}=U_{n-1}+U_{n-2}+U_{n-3}+U_{n-4},\text{ \ \ \ \ }%
U_{0}=U_{1}=0,U_{2}=U_{3}=1.
\end{equation*}%
The numbers $U_{n}$ can be expressed using Binet's formula%
\begin{equation*}
U_{n}=\frac{\alpha ^{n}}{(\alpha -\beta )(\alpha -\gamma )(\alpha -\delta )}+%
\frac{\beta ^{n}}{(\beta -\alpha )(\beta -\gamma )(\beta -\delta )}+\frac{%
\gamma ^{n}}{(\gamma -\alpha )(\gamma -\beta )(\gamma -\delta )}+\frac{%
\delta ^{n}}{(\delta -\alpha )(\delta -\beta )(\delta -\gamma )}.
\end{equation*}%
We define the square matrix $B$ of order $4$ as:%
\begin{equation*}
B=\left( 
\begin{array}{cccc}
1 & 1 & 1 & 1 \\ 
1 & 0 & 0 & 0 \\ 
0 & 1 & 0 & 0 \\ 
0 & 0 & 1 & 0%
\end{array}%
\right)
\end{equation*}%
such that $\det B=-1.$

Induction proof may be used to establish 
\begin{equation}
B^{n}=\left( 
\begin{array}{cccc}
U_{n+2} & U_{n+1}+U_{n}+U_{n-1} & U_{n+1}+U_{n} & U_{n+1} \\ 
U_{n+1} & U_{n}+U_{n-1}+U_{n-2} & U_{n}+U_{n-1} & U_{n} \\ 
U_{n} & U_{n-1}+U_{n-2}+U_{n-3} & U_{n-1}+U_{n-2} & U_{n-1} \\ 
U_{n-1} & U_{n-2}+U_{n-3}+U_{n-4} & U_{n-2}+U_{n-3} & U_{n-2}%
\end{array}%
\right) .  \label{equation:gbhasmnuopdcxsz}
\end{equation}%
Matrix formulation of $M_{n}$ and $R_{n}$ can be given as%
\begin{equation}
\left( 
\begin{array}{c}
M_{n+3} \\ 
M_{n+2} \\ 
M_{n+1} \\ 
M_{n}%
\end{array}%
\right) =\left( 
\begin{array}{cccc}
1 & 1 & 1 & 1 \\ 
1 & 0 & 0 & 0 \\ 
0 & 1 & 0 & 0 \\ 
0 & 0 & 1 & 0%
\end{array}%
\right) ^{n}\left( 
\begin{array}{c}
M_{3} \\ 
M_{2} \\ 
M_{1} \\ 
M_{0}%
\end{array}%
\right)  \label{equat:nmouyvbcfxdsz}
\end{equation}%
and%
\begin{equation}
\left( 
\begin{array}{c}
R_{n+3} \\ 
R_{n+2} \\ 
R_{n+1} \\ 
R_{n}%
\end{array}%
\right) =\left( 
\begin{array}{cccc}
1 & 1 & 1 & 1 \\ 
1 & 0 & 0 & 0 \\ 
0 & 1 & 0 & 0 \\ 
0 & 0 & 1 & 0%
\end{array}%
\right) ^{n}\left( 
\begin{array}{c}
R_{3} \\ 
R_{2} \\ 
R_{1} \\ 
R_{0}%
\end{array}%
\right) .  \label{equatio:mnbvuhnbtdfx}
\end{equation}%
Induction proofs may be used to establish the matrix formulations $M_{n}$
and $R_{n}$.

Now we define the matrices $B_{M}$ and $B_{R}$ as 
\begin{equation*}
B_{M}=\left( 
\begin{array}{cccc}
\widehat{M}_{5} & \widehat{M}_{4}+\widehat{M}_{3}+\widehat{M}_{2} & \widehat{%
M}_{4}+\widehat{M}_{3} & \widehat{M}_{4} \\ 
\widehat{M}_{4} & \widehat{M}_{3}+\widehat{M}_{2}+\widehat{M}_{1} & \widehat{%
M}_{3}+\widehat{M}_{2} & \widehat{M}_{3} \\ 
\widehat{M}_{3} & \widehat{M}_{2}+\widehat{M}_{1}+\widehat{M}_{0} & \widehat{%
M}_{2}+\widehat{M}_{1} & \widehat{M}_{2} \\ 
\widehat{M}_{2} & \widehat{M}_{1}+\widehat{M}_{0}+\widehat{M}_{-1} & 
\widehat{M}_{1}+\widehat{M}_{0} & \widehat{M}_{1}%
\end{array}%
\right) \text{\ and }B_{R}=\left( 
\begin{array}{cccc}
\widehat{R}_{5} & \widehat{R}_{4}+\widehat{R}_{3}+\widehat{R}_{2} & \widehat{%
R}_{4}+\widehat{R}_{3} & \widehat{R}_{4} \\ 
\widehat{R}_{4} & \widehat{R}_{3}+\widehat{R}_{2}+\widehat{R}_{1} & \widehat{%
R}_{3}+\widehat{R}_{2} & \widehat{R}_{3} \\ 
\widehat{R}_{3} & \widehat{R}_{2}+\widehat{R}_{1}+\widehat{R}_{0} & \widehat{%
R}_{2}+\widehat{R}_{1} & \widehat{R}_{2} \\ 
\widehat{R}_{2} & \widehat{R}_{1}+\widehat{R}_{0}+\widehat{R}_{-1} & 
\widehat{R}_{1}+\widehat{R}_{0} & \widehat{R}_{1}%
\end{array}%
\right) .
\end{equation*}%
These matrices\ $B_{M}$ and $B_{R}$ can be called Tetranacci quaternion
matrix and Tetranacci-Lucas quaternion matrix, respectively.

\begin{theorem}
For $n\geq 0,$ the followings are valid:

\begin{description}
\item[(a)] 
\begin{equation}
B_{M}\left( 
\begin{array}{cccc}
1 & 1 & 1 & 1 \\ 
1 & 0 & 0 & 0 \\ 
0 & 1 & 0 & 0 \\ 
0 & 0 & 1 & 0%
\end{array}%
\right) ^{n}=\left( 
\begin{array}{cccc}
\widehat{M}_{n+5} & \widehat{M}_{n+4}+\widehat{M}_{n+3}+\widehat{M}_{n+2} & 
\widehat{M}_{n+4}+\widehat{M}_{n+3} & \widehat{M}_{n+4} \\ 
\widehat{M}_{n+4} & \widehat{M}_{n+3}+\widehat{M}_{n+2}+\widehat{M}_{n+1} & 
\widehat{M}_{n+3}+\widehat{M}_{n+2} & \widehat{M}_{n+3} \\ 
\widehat{M}_{n+3} & \widehat{M}_{n+2}+\widehat{M}_{n+1}+\widehat{M}_{n} & 
\widehat{M}_{n+2}+\widehat{M}_{n+1} & \widehat{M}_{n+2} \\ 
\widehat{M}_{n+2} & \widehat{M}_{n+1}+\widehat{M}_{n}+\widehat{M}_{n-1} & 
\widehat{M}_{n+1}+\widehat{M}_{n} & \widehat{M}_{n+1}%
\end{array}%
\right) ,  \label{equation:nbhbxczdradweq}
\end{equation}

\item[(b)] 
\begin{equation}
B_{R}\left( 
\begin{array}{cccc}
1 & 1 & 1 & 1 \\ 
1 & 0 & 0 & 0 \\ 
0 & 1 & 0 & 0 \\ 
0 & 0 & 1 & 0%
\end{array}%
\right) ^{n}=\left( 
\begin{array}{cccc}
\widehat{R}_{n+5} & \widehat{R}_{n+4}+\widehat{R}_{n+3}+\widehat{R}_{n+2} & 
\widehat{R}_{n+4}+\widehat{R}_{n+3} & \widehat{R}_{n+4} \\ 
\widehat{R}_{n+4} & \widehat{R}_{n+3}+\widehat{R}_{n+2}+\widehat{R}_{n+1} & 
\widehat{R}_{n+3}+\widehat{R}_{n+2} & \widehat{R}_{n+3} \\ 
\widehat{R}_{n+3} & \widehat{R}_{n+2}+\widehat{R}_{n+1}+\widehat{R}_{n} & 
\widehat{R}_{n+2}+\widehat{R}_{n+1} & \widehat{R}_{n+2} \\ 
\widehat{R}_{n+2} & \widehat{R}_{n+1}+\widehat{R}_{n}+\widehat{R}_{n-1} & 
\widehat{R}_{n+1}+\widehat{R}_{n} & \widehat{R}_{n+1}%
\end{array}%
\right) .  \label{equa:tyghbcvsdae}
\end{equation}
\end{description}
\end{theorem}

\textit{Proof}. We prove (a) by mathematical induction on $n.$ If $n=0,$
then the result is clear. Now, we assume it is true for $n=k,$ that is%
\begin{equation*}
B_{M}B^{k}=\left( 
\begin{array}{cccc}
\widehat{M}_{k+5} & \widehat{M}_{k+4}+\widehat{M}_{k+3}+\widehat{M}_{k+2} & 
\widehat{M}_{k+4}+\widehat{M}_{k+3} & \widehat{M}_{k+4} \\ 
\widehat{M}_{k+4} & \widehat{M}_{k+3}+\widehat{M}_{k+2}+\widehat{M}_{k+1} & 
\widehat{M}_{k+3}+\widehat{M}_{k+2} & \widehat{M}_{k+3} \\ 
\widehat{M}_{k+3} & \widehat{M}_{k+2}+\widehat{M}_{k+1}+\widehat{M}_{k} & 
\widehat{M}_{k+2}+\widehat{M}_{k+1} & \widehat{M}_{k+2} \\ 
\widehat{M}_{k+2} & \widehat{M}_{k+1}+\widehat{M}_{k}+\widehat{M}_{k-1} & 
\widehat{M}_{k+1}+\widehat{M}_{k} & \widehat{M}_{k+1}%
\end{array}%
\right) .
\end{equation*}%
If we use (\ref{equation:abgsvbnuytsaerds}), then we have $\widehat{M}_{k+4}=%
\widehat{M}_{k+3}+\widehat{M}_{k+2}+\widehat{M}_{k+1}+\widehat{M}_{k}.$
Then, by induction hypothesis, we obtain%
\begin{eqnarray*}
B_{M}B^{k+1} &=&(B_{M}B^{k})B \\
&=&\left( 
\begin{array}{cccc}
\widehat{M}_{k+5} & \widehat{M}_{k+4}+\widehat{M}_{k+3}+\widehat{M}_{k+2} & 
\widehat{M}_{k+4}+\widehat{M}_{k+3} & \widehat{M}_{k+4} \\ 
\widehat{M}_{k+4} & \widehat{M}_{k+3}+\widehat{M}_{k+2}+\widehat{M}_{k+1} & 
\widehat{M}_{k+3}+\widehat{M}_{k+2} & \widehat{M}_{k+3} \\ 
\widehat{M}_{k+3} & \widehat{M}_{k+2}+\widehat{M}_{k+1}+\widehat{M}_{k} & 
\widehat{M}_{k+2}+\widehat{M}_{k+1} & \widehat{M}_{k+2} \\ 
\widehat{M}_{k+2} & \widehat{M}_{k+1}+\widehat{M}_{k}+\widehat{M}_{k-1} & 
\widehat{M}_{k+1}+\widehat{M}_{k} & \widehat{M}_{k+1}%
\end{array}%
\right) \left( 
\begin{array}{cccc}
1 & 1 & 1 & 1 \\ 
1 & 0 & 0 & 0 \\ 
0 & 1 & 0 & 0 \\ 
0 & 0 & 1 & 0%
\end{array}%
\right)  \\
&=&\left( 
\begin{array}{cccc}
\widehat{M}_{k+5}+\widehat{M}_{k+4}+\widehat{M}_{k+3}+\widehat{M}_{k+2} & 
\widehat{M}_{k+5}+\widehat{M}_{k+4}+\widehat{M}_{k+3} & \widehat{M}_{k+5}+%
\widehat{M}_{k+4} & \widehat{M}_{k+5} \\ 
\widehat{M}_{k+4}+\widehat{M}_{k+3}+\widehat{M}_{k+2}+\widehat{M}_{k+1} & 
\widehat{M}_{k+4}+\widehat{M}_{k+3}+\widehat{M}_{k+2} & \widehat{M}_{k+4}+%
\widehat{M}_{k+3} & \widehat{M}_{k+4} \\ 
\widehat{M}_{k+3}+\widehat{M}_{k+2}+\widehat{M}_{k+1}+\widehat{M}_{k} & 
\widehat{M}_{k+3}+\widehat{M}_{k+2}+\widehat{M}_{k+1} & \widehat{M}_{k+3}+%
\widehat{M}_{k+2} & \widehat{M}_{k+3} \\ 
\widehat{M}_{k+2}+\widehat{M}_{k+1}+\widehat{M}_{k}+\widehat{M}_{k-1} & 
\widehat{M}_{k+2}+\widehat{M}_{k+1}+\widehat{M}_{k} & \widehat{M}_{k+2}+%
\widehat{M}_{k+1} & \widehat{M}_{k+2}%
\end{array}%
\right)  \\
&=&\left( 
\begin{array}{cccc}
\widehat{M}_{k+6} & \widehat{M}_{k+5}+\widehat{M}_{k+4}+\widehat{M}_{k+3} & 
\widehat{M}_{k+5}+\widehat{M}_{k+4} & \widehat{M}_{k+5} \\ 
\widehat{M}_{k+5} & \widehat{M}_{k+4}+\widehat{M}_{k+3}+\widehat{M}_{k+2} & 
\widehat{M}_{k+4}+\widehat{M}_{k+3} & \widehat{M}_{k+4} \\ 
\widehat{M}_{k+4} & \widehat{M}_{k+3}+\widehat{M}_{k+2}+\widehat{M}_{k+1} & 
\widehat{M}_{k+3}+\widehat{M}_{k+2} & \widehat{M}_{k+3} \\ 
\widehat{M}_{k+3} & \widehat{M}_{k+2}+\widehat{M}_{k+1}+\widehat{M}_{k} & 
\widehat{M}_{k+2}+\widehat{M}_{k+1} & \widehat{M}_{k+2}%
\end{array}%
\right) .
\end{eqnarray*}%
Thus, (\ref{equation:nbhbxczdradweq}) holds for all non-negative integers $n.
$

(\ref{equa:tyghbcvsdae}) can be similarly proved . 
\endproof%

\begin{corollary}
For $n\geq 0,$ the followings hold:

\begin{description}
\item[(a)] $\widehat{M}_{n+3}=\widehat{M}_{3}U_{n+2}+(\widehat{M}_{2}+%
\widehat{M}_{1}+\widehat{M}_{0})U_{n+1}+(\widehat{M}_{1}+\widehat{M}%
_{2})U_{n}+\widehat{M}_{2}U_{n-1}$

\item[(b)] $\widehat{R}_{n+3}=\widehat{R}_{3}U_{n+2}+(\widehat{R}_{2}+%
\widehat{R}_{1}+\widehat{R}_{0})U_{n+1}+(\widehat{R}_{1}+\widehat{R}%
_{2})U_{n}+\widehat{R}_{2}U_{n-1}$
\end{description}
\end{corollary}

\textit{Proof}. The proof of (a) can be seen by the coefficient of the
matrix $B_{M}$ and (\ref{equation:gbhasmnuopdcxsz}). The proof of (b) can be
seen by the coefficient\ of the matrix $B_{R}$ and (\ref%
{equation:gbhasmnuopdcxsz}). 
\endproof%

\end{document}